\documentclass[a4paper]{amsart}
\usepackage{amsmath,amsthm,amssymb}
\usepackage[pdftex]{graphicx}
\usepackage{mathrsfs}
\usepackage[dvipdfmx, usenames]{color}
\usepackage[all]{xy}
\usepackage{comment}

\newtheorem{thm}{Theorem}[section]

\newtheorem{prop}[thm]{Proposition}

\theoremstyle{definition}

\newtheorem{defn}[thm]{Definition}

\theoremstyle{remark}

\newtheorem{rem}[thm]{Remark}        

\numberwithin{equation}{section}

\def\XXint#1#2#3{{\setbox0=\hbox{$#1{#2#3}{\int}$}
\vcenter{\hbox{$#2#3$}}\kern-.5\wd0}}

\newcommand{\R}{\mathbb{R}}

\newcommand{\m}{\mathfrak{m}}

\newcommand{\Cpl}{\mathsf{Cpl}}

\newcommand{\supp}{\mathsf{supp}\,}

\newcommand{\LIP}{{\bf \mathrm{LIP}}}

\newcommand{\lip}{\mathrm{lip}}
\newcommand{\Ch}{\mathsf{Ch}}

\newcommand{\RCD}{\mathsf{RCD}}
\let\CD\relax
\newcommand{\CD}{\mathsf{CD}}
\newcommand{\ul}{\underline}

\newcommand{\la}{\left\langle}
\newcommand{\ra}{\right\rangle}

\newcommand{\lv}{\left\vert}
\newcommand{\rv}{\right\vert}
\newcommand{\lV}{\left\Vert}
\newcommand{\rV}{\right\Vert}

\newcommand{\dime}{\mathsf{dim}_{\mathsf{ess}}}
\newcommand{\Ent}{\mathsf{Ent}}
\newcommand{\Tan}{\mathsf{Tan}}

\newcommand{\calD}{\mathcal{D}}

\newcommand{\calH}{\mathcal{H}}

\newcommand{\calL}{\mathcal{L}}

\newcommand{\calP}{\mathcal{P}}

\newcommand{\calR}{\mathcal{R}}
\newcommand{\calS}{\mathcal{S}}

\newcommand{\bbN}{\mathbb{N}}

\usepackage[abbrev,nobysame]{amsrefs}

\makeatletter
\def\@makefnmark{%
\leavevmode
\raise.9ex\hbox{\check@mathfonts
\fontsize\sf@size\z@\normalfont%
\@thefnmark}%
}
\makeatother

\title{One dimensional $\mathsf{RCD}$ spaces always satisfy the regular Weyl's law}
\author{Akemi Iwahashi, Yu Kitabeppu, Akari Yonekura}
\address{Kumamoto University}
\email[Akemi Iwahashi]{211d8001@st.kumamoto-u.ac.jp}
\email[Yu Kitabeppu]{ybeppu@kumamoto-u.ac.jp}
\email[Akari Yonekura]{216d8018@st.kumamoto-u.ac.jp}

\begin{document}
\maketitle
 \begin{abstract}
  Ambrosio, Honda, and Tewodrose proved that the regular Weyl's law is equivalent to a mild condition related to the infinitesimal behavior of the measure of balls in compact finite dimensional $\RCD$ spaces. Though that condition is seemed to always hold for any such spaces, however, Dai, Honda, Pan, and Wei recently show that for any integer $n$ at least 2, there exists a compact $\RCD$ space of $n$ dimension fails to satisfy the regular Weyl's law. In this short article we prove that one dimensional $\RCD$ spaces always satisfy the regular Weyl's law. 
 \end{abstract}
%
%
\section{Introduction}
In 1911, H. Weyl proved an asymptotic formula of the distribution of the eigenvalues of the Dirichlet Laplacian in the bounded domain $\Omega\subset\R^n$(\cite{W}); precisely, let $\{\lambda_i\}_i$ be a sequence of all eigenvalues of the (minus) Dirichlet Laplacian $-\Delta$, and let $N(\lambda)$ be the counting function of the eigenvalues, namely $N(\lambda):=\#\{i\;;\;\lambda_i\leq\lambda\}$. It is known that $\lambda_i$ are positive numbers and $\lambda_i\rightarrow \infty$. Then it holds that 
\begin{align}
 \lim_{\lambda\rightarrow\infty}\frac{N(\lambda)}{\lambda^{n/2}}=\frac{\omega_n}{(2\pi)^n}\calH^n(\Omega),\notag
\end{align} 
where $\omega_n$ is the volume of the unit ball in $\R^n$ and $\calH^n$ is the $n$-dimensional Hausdorff measure. This result has been generalized to many other situations(see for instance \cite{I100} for a brief history and the generalization of that theorems). Here we focus on non-smooth setting. In \cite{AHTweyl}, Ambrosio, Honda, and Tewodrose proved the following theorem. 
\begin{thm}
For a compact finite dimensional $\RCD$ space $(X,d,\m)$ of the essential dimension $n$(see Definition \ref{def:essdim} for the definition), the Weyl's law 
\begin{align}
 \lim_{\lambda\rightarrow\infty}\frac{N(\lambda)}{\lambda^{n/2}}=\frac{\omega_n}{(2\pi)^n}\calH^n(X)\label{eq:RegWeyl}
\end{align}
holds if and only if 
\begin{align}
 \lim_{r\rightarrow +0}\int_X\frac{r^n}{\m(B_r(x))}\,\m(dx)=\int_X\lim_{r\rightarrow+0}\frac{r^n}{\m(B_r(x))}\,\m(dx).\label{eq:commuteintlim}
\end{align}
\end{thm}
We give some remarks. A finite dimensional $\RCD$ space $(X,d,\m)$ is a metric measure space with Ricci curvature bounded from below and dimension bounded from above in the synthetic sense(see Definition \ref{defn:RCD} for precise definition).  
Under (\ref{eq:commuteintlim}), though the reference measure $\m$ is not $\calH^n$ in general, the $n$-dimensional Hausdorff measure appears in the right-hand side in (\ref{eq:RegWeyl}).    

Since it is known that $\lim_{r\rightarrow+0}r^n/\m(B_r(x))\in(0,\infty)$ $\m$-a.e.(\cite{AHTweyl}), (\ref{eq:commuteintlim}) seems to be held for \emph{all} $\RCD$ spaces. However, Dai, Honda, Pan, and Wei prove the following. 
 \begin{thm}[\cite{DHPW}]
  For any $\beta\in(2,\infty)$, there exists a compact $\RCD(-1,N_{\beta})$ space $(X,d,\m)$ of the essential dimension 2 such that 
  \begin{align}
   \lim_{\lambda\rightarrow \infty}\frac{N(\lambda)}{\lambda^{\beta/2}}=\frac{c}{\Gamma(\beta+1)}\calH^{\beta}(\calS)\in(0,\infty),\notag
  \end{align}
  where $\calS\subset X$ is the singular set and $c$ is a canonical constant.   
 \end{thm}
See Theorem \ref{thm:BS} for the definition of the singular set and \cite{DHPW} the explicit definition of the canonical constant $c$. One notably fact is that the limit of the counting function is different from the ordinary ones. Moreover, they show a more pathological example. 
 \begin{thm}[\cite{DHPW}]
  There exists a compact $\RCD(-1,10)$ space such that 
  \begin{align}
   \lim_{\lambda\rightarrow \infty}\frac{N(\lambda)}{\lambda\log\lambda}=\frac{1}{4\pi}.\notag
  \end{align}
 \end{thm}
In this example, the asymptotic behavior of the counting function is not even polynomial growth. 

Motivated by these theorem, (\ref{eq:RegWeyl}) is called the \emph{regular Weyl's law}. So a natural question is raised: 
\medskip
\par $\ul{Q}$: Does the regular Weyl's law hold for one-dimensional $\RCD$ spaces? 
\medskip 
\par The following is the main theorem. 
 \begin{thm}
  Let $(X,d,\m)$ be a non-trivial compact $\RCD(K,N)$ space for $K\in\R$, $N\in(1,\infty)$. Then the followings are all equivalent; 
  \begin{enumerate}
   \item $\dime (X,d,\m)=1$, 
   \item $X=[0,\ell]$ or $S^1(r)$ for $\ell,r>0$, 
   \item $N(\lambda)\sim \lambda^{1/2}$, 
   \item $\displaystyle{\lim_{\lambda\downarrow 0}\frac{N(\lambda)}{\lambda^{1/2}}=\frac{\omega_1}{2\pi}\calH^1(X)=\frac{1}{\pi}\calH^1(X)}$, 
   \item $\displaystyle{\lim_{\lambda\rightarrow\infty}\frac{N(\lambda)}{\lambda^{(1+\alpha)/2}}}=0$ for $0<\alpha\leq 1$. 
  \end{enumerate}
 \end{thm}
This theorem is a complete answer to the question. And on the contrary for $\RCD$ spaces of dimension at least 2, the regular Weyl's law always holds for \emph{ALL} one-dimensional $\RCD$ spaces.

%
%
\section{Preliminaries}
Let $(Y,d_Y)$ be a metric space. We denote the set of all continuous functions on $Y$ with bounded support by $C_{bs}(Y)$. We call a function $g:Y\rightarrow \R$ \emph{an $L$-Lipschitz function} for $L>0$ if $\lv g(y_0)-g(y_1)\rv\leq Ld_Y(y_0,y_1)$ holds for any $y_0,y_1\in Y$. We denote the set of all Lipschitz functions in $Y$ by $\LIP(Y)$. For a generic function $f:Y\rightarrow \R$, its \emph{local Lipschitz constant} at $y\in Y$ is defined by 
\begin{align}
 \lip f(y):=\begin{cases}
  \limsup_{z\rightarrow y}\frac{\lv f(z)-f(y)\rv}{d_Y(z,y)}&\text{if $y$ is not isolated},\\
  0&\text{otherwise}.
 \end{cases}\notag
\end{align}
We denote the set of all Borel probability measures on $Y$ by $\calP(Y)$. We define 
\begin{align}
 \calP_2(Y):=\left\{\mu\in\calP(Y)\;;\;\int_Yd_Y^2(~^{\exists}o,y)\,\mu(dy)<\infty\right\}. \notag
\end{align} 
For any two Borel probability measures $\mu,\nu\in\calP(Y)$, \emph{the coupling} $\xi\in\calP(Y\times Y)$ between them is defined as 
\begin{align}
 \begin{cases}
  \xi(A\times Y)=\mu(A)\\
  \xi(Y\times A)=\nu(A)&\text{for any Borel subset }A\subset Y. 
 \end{cases}\notag
\end{align}
The set of all couplings between $\mu$ and $\nu$ is denoted by $\Cpl(\mu,\nu)$, which is not empty since $\mu\otimes\nu\in\Cpl(\mu,\nu)$. 
\begin{defn}[$L^2$-Wasserstein space]
 For given $\mu,\nu\in\calP_2(Y)$, \emph{the $L^2$-Wasserstein distance} between them, $W_2(\mu,\nu)$, is defined by 
 \begin{align}
  W_2(\mu,\nu):=\inf\left\{\lV d\rV_{L^2(\xi)}\;;\;\xi\in\Cpl(\mu,\nu)\right\}.\notag
 \end{align} 
 It is known that $W_2$ is a metric on $\calP_2(Y)$. The metric space $(\calP_2(Y),W_2)$ is called \emph{the $L^2$-Wasserstein space}. 
\end{defn}
It is known that $(\calP_2(Y),W_2)$ is complete separable if and only if so is $(Y,d_Y)$.  
\subsection{Convex functions on geodesic spaces}
Let $(Y,d_Y)$ be a metric space. We call $(Y,d_Y)$ \emph{a geodesic space} if for any two points $y_0,y_1\in Y$, there exists a continuous curve $\gamma:[0,1]\rightarrow Y$ connecting them such that 
\begin{align}
 d_Y(\gamma_s,\gamma_t)=\lv s-t\rv d_Y(y_0,y_1)\notag
\end{align}
holds for any $s,t\in [0,1]$. 
For given $K\in\R$, $N\in(1,\infty)$, we define the \emph{distortion coefficients} $\sigma^{(t)}_{K,N}$ for $t\in[0,1]$ by 
\begin{align}
 \sigma^{(t)}_{K,N}(\theta):=\begin{cases}
 \infty&\text{if }K\theta^2\geq N\pi^2,\\
 \frac{\sin(t\theta\sqrt{K/N})}{\sin(\theta)\sqrt{K/N}}&\text{if }0<K\theta^2<N\pi^2,\\
 t&\text{if }K\theta^2=0,\\
 \frac{\sinh(t\theta\sqrt{-K/N})}{\sinh(\theta\sqrt{-K/N})}&\text{if }K\theta^2<0.
 \end{cases}\notag
\end{align}
By using the distortion coefficients, we define the convexity of functions. For a given function $g:Y\rightarrow \R\cup\{\infty\}$, $\calD(g):=\{y\in Y\;;\; g(y)<\infty\}$. 
 \begin{defn}[$(K,N)$-convex]
  We say a function $f:Y\rightarrow \R\cup\{\infty\}$, \emph{$(K,N)$-convex} if for any $y_0,y_1\in Y$, there exists a geodesic $\gamma:[0,1]\rightarrow Y$ connecting them such that 
  \begin{align}
   &\exp\left(-\frac{1}{N}f(\gamma_t)\right)\notag\\
   &\geq \sigma^{(1-t)}_{K,N}(d_Y(y_0,y_1))\exp\left(-\frac{1}{N}f(y_0)\right)+\sigma^{(t)}_{K,N}(d_Y(y_0,y_1))\exp\left(-\frac{1}{N}f(y_1)\right)\notag
  \end{align}
  for any $t\in[0,1]$. 
 \end{defn}
One can prove that each $(K,N)$-convex function is locally Lipschitz in the interior of the geodesic. Hence $(K,N)$-convex functions on the intervals are differentiable almost everywhere by Rademacher's theorem.  
The following results, that is a key proposition in this short article, are found in \cite{CMnew} (cf. \cite{CMonge}). 
 \begin{prop}[\cite{CMnew}] 
  Let $I:=(a,b)$ be an interval for $a<b$ and $f$ a $(K,N-1)$-convex function for $K<0$ on $I$. Then any $x_0<x_1\in I$, 
  \begin{align}
   \left(\frac{\sinh((b-x_1)\sqrt{-K/N-1})}{\sinh((b-x_0)\sqrt{-K/N-1})}\right)^{N-1}\leq \frac{e^{-f(x_1)}}{e^{-f(x_0)}}\leq\left(\frac{\sinh((x_1-a)\sqrt{-K/N-1})}{\sinh((x_0-a)\sqrt{-K/N-1})}\right)^{N-1}.\label{eq:hratiioest}
  \end{align}
 \end{prop}
\subsection{$\RCD$ spaces}
We call a triplet $(X,d,\m)$ metric measure space if $(X,d)$ is a complete separable metric space and $\m$ is a locally finite Borel measure on $X$. Two metric measure spaces $(X,d_X,\m_X)$ and $(Y,d_Y,\m_Y)$ are isomorphic to each other if there exists an isometry $f:\supp\m_X\rightarrow \supp\m_Y$ such that $f_*\m_X=\m_Y$. In this case we just denote $(X,d_X,\m_X)=(Y,d_Y,\m_Y)$. 
For a metric measure space $(X,d,\m)$, $(X,\m)$ is a $\sigma$-finite measure space. Hence the Radon-Nikodym theorem holds, that is, if $\mu\in\calP(X)$ is absolutely continuous with respect to $\m$, denoted by $\mu\ll\m$, there exists an $L^1(\m)$ function $\rho$ such that $\mu=\rho\m$. From now on, we always assume that there exists a constant $C>0$ such that 
\begin{align}
 \int_Xe^{-Cd^2(x_0,x)}\,\m(dx)<\infty\label{eq:basicvolrate}
\end{align}
for a point $x_0\in X$. 
\emph{The relative entropy} functional $\Ent_{\m}:\calP(X)\rightarrow \R\cup\{\pm\infty\}$ is defined by 
\begin{align}
 \Ent_{\m}(\mu):=\begin{cases}
  \int_{\{\rho>0\}}\rho\log\rho\,d\m&\text{if }\mu=\rho\m\ll\m\\
  \infty&\text{otherwise}. 
 \end{cases}\notag
\end{align}
By (\ref{eq:basicvolrate}), $\Ent_{\m}(\mu)>-\infty$ for any $\mu\in\calP_2(X)$. 
 \begin{defn}[$\CD_e(K,N)$ space,\cite{EKS}]
  Let $K\in\R$, $N\in(1,\infty)$. A metric measure space $(X,d,\m)$ with (\ref{eq:basicvolrate}) is called \emph{a $\CD_e(K,N)$ space} if $\Ent_{\m}$ is $(K,N)$-convex. 
 \end{defn}
For $f\in L^2(\m)$, we define \emph{the Cheeger energy} of $f$, $\Ch(f)$, by 
\begin{align}
 \Ch(f):=\frac{1}{2}\inf\left\{\liminf_{n\rightarrow\infty}\int_X(\lip f_n)^2\,d\m\;;\;f_n\xrightarrow{L^2(\m)} f,\,f_n\in\LIP(X)\right\}.\notag
\end{align} 
It is known that for any $f\in \calD(\Ch)$, there exists an $L^2(\m)$-function $\lv df\rv$ such that $2\Ch(f)=\int_X\lv df\rv^2\,d\m$. We define a Banach space $W^{1,2}(X,d,\m):=L^2(\m)\cap \calD(\Ch)$ equipped with the norm 
\begin{align}
 \lV f\rV^2_{1,2}:=\lV f\rV^2_{L^2(\m)}+\lV \lv df\rv\rV^2_{L^2(\m)},\notag
\end{align} 
called \emph{the Sobolev space} over $(X,d,\m)$. 
 \begin{defn}[Infinitesimal Hilbertianity, \cite{Gdiffstr}]
  Let $(X,d,\m)$ be a metric measure space. We say that $(X,d,\m)$ is \emph{infinitesimally Hilbertian} if $W^{1,2}(X,d,\m)$ is a Hilbert space.  
 \end{defn}
By the infinitesimal Hilbertianity, we define the inner product of the differentials for $f,g\in\calD(\Ch)$ by 
\begin{align}
 \int_X\la  df,dg\ra\,d\m:=\frac{1}{2}\left(\Ch(f+g)-\Ch(f-g)\right).\notag
\end{align}
 \begin{defn}[$\RCD$ space]\label{defn:RCD}
  Let $(X,d,\m)$ be a metric measure space with the volume growth condition (\ref{eq:basicvolrate}), and $K\in\R$, $N\in(1,\infty)$. We call $(X,d,\m)$ an $\RCD(K,N)$ space if the following two conditions are hold; 
  \begin{enumerate}
   \item $(X,d,\m)$ is a $\CD_e(K,N)$ space, 
   \item $(X,d,\m)$ is infinitesimally Hilbertian. 
  \end{enumerate}
 \end{defn}
The $n$-dimensional Riemannian manifolds with $\mathrm{Ric}\geq K$ are $\RCD(K,N)$ spaces if $n\leq N$, and Ricci limit spaces, measured Gromov-Hausdorff limit of manifolds with bounded dimension from above and bounded Ricci curvature from below, are also $\RCD$ spaces. See Theorem \ref{thm:RCDstable}.  
 \begin{rem}[Historical remarks]
 The synthetic notion of Ricci curvature bound for metric measure spaces is first introduced by Sturm \cite{Stmms1,Stmms2} and Lott, Villani \cite{LV} independently(\cite{LV} for $\CD(K,\infty)$ and $\CD(0,N)$ for finite $N$, \cite{Stmms1} for $\CD(K,\infty)$, \cite{Stmms2} for $\CD(K,N)$ for finite $N$). In order to get the tensorial and localization property, Bacher and Sturm introduced the reduced curvature-dimension condition $\CD^*(K,N)$ \cite{BS}. Non-Riemannian Finsler manifolds can be $\CD$ spaces. To get rid of such a class of spaces, Gigli defined the infinitesimal Hilbertianity \cite{Gdiffstr} and Ambrosio, Gigli, Savar\'e defined the Riemannian curvature-dimension condition $\RCD(K,\infty)$ for compact metric measure spaces \cite{AGSRiem}, afterwards, Ambrosio, Gigli, Mondino, Rajala defined the same notion for $\sigma$-finite cases \cite{AGMR}. For finite $N$, Erbar, Kuwada, Sturm \cite{EKS} and Ambrosio, Mondino, Savar\'e \cite{AMSnon} defined $\RCD^*(K,N)$ space independently. Under essentially non-branching assumption(see \cite{RaS} for the definition), $\CD$ condition and $\CD^*$ condition are equivalent to each other\cite{CMil}. The metric measure spaces that satisfy the conditions in Definition \ref{defn:RCD} was called $\RCD^*(K,N)$ spaces. However, by \cite{CMil}, we call these spaces $\RCD(K,N)$ spaces now. 
 
 In order to emphasize $N$ being finite, we say that $\RCD(K,N)$ space is \emph{finite dimensional}. 
 \end{rem}
 \begin{rem}\label{rem:domnmul}
  The $\RCD(K,N)$ condition is the synthetic notion of lower bound of Ricci curvature($\geq K$), and upper bound of dimension($\leq N$). Actually, $\RCD(K,N)$ space $(X,d,\m)$ also satisfies $\RCD(K',N')$ condition for $K'\leq K$ and $N'\geq N$.  
  \smallskip
  \par For $\RCD(K,N)$ space $(X,d,\m)$, $(X,ad,b\m)$ is a $\RCD(a^{-2}K,N)$ space for $a,b>0$.   
 \end{rem}

Let $(X,d,\m)$ be an $\RCD(K,N)$ space. 
\begin{defn}[Laplacian]
 A function $f\in \calD(\Ch)$ belongs to $\calD(\Delta)$ if there exists an $L^2(\m)$-function $h$ such that 
 \begin{align}
  \int_X\la df,dg\ra\,d\m=-\int_Xhg\,d\m\notag
 \end{align}
 holds for any $g\in\calD(\Ch)$. In this case, the Laplacian of $f$ is denoted by $\Delta f:=h$. 
\end{defn}
The Laplacian is a densely defined nonpositively definite self-adjoint operator in $L^2(\m)$. When $(X,d,\m)$ is a compact finite dimensional $\RCD$ space, all the spectrum of $-\Delta$ are eigenvalues, and the multiplicity of zero spectrum is 1 since the resolvent operators are compact. We denote the non-zero eigenvalues by $0<\lambda_1\leq \lambda_2\leq\cdots\rightarrow\infty$ with multiplicity.  
\subsection{Infinitesimal structure on $\RCD$ spaces}
Let $\{(X_i,d_i,\m_i,x_i^0)\}_{i\in\bbN\cup\{\infty\}}$ be a family of pointed metric measure spaces. We say that $(X_i,d_i,\m_i,x^0_i)$ converges to $(X_{\infty},d_{\infty},\m_{\infty},x^0_{\infty})$ in \emph{the measured Gromov-Hausdorff} sense (mGH for short) if the following conditions are satisfied; there exist sequences of positive numbers $\varepsilon_i\downarrow 0$, $R_i\uparrow\infty$, and of Borel maps $\varphi_i:B_{R_i}(x_i^0)\rightarrow X_{\infty}$
\begin{enumerate}
 \item $\lv d_i(x,y)-d_{\infty}(\varphi_i(x),\varphi_i(y))\rv<\varepsilon_i$ for any $i\in\bbN$, any $x,y\in B_{R_i}(x_i^0)$. And $B_{R_i-\varepsilon_i}(\varphi_i(x_i))\subset B_{\varepsilon_i}(\varphi_i(B_{R_i}(x_i^0)))$, 
 \item $\varphi_i(x_i^0)\rightarrow x^0_{\infty}$, 
 \item $(\varphi_i)_{*}\m_i\rightharpoonup \m_{\infty}$ in duality with $C_{bs}(X_{\infty})$. 
\end{enumerate}
In this case we denote by $(X_i,d_i,\m_i,x_i^0)\xrightarrow{mGH}(X_{\infty},d_{\infty},\m_{\infty},x_{\infty}^0)$. $\RCD$ condition is stable under mGH convergence and a set of pointed metric measure spaces with the same $\RCD$ condition is sequentially precompact with respect to mGH convergence(see \cite{GMS}).  
 \begin{thm}[\cite{GMS,EKS}]\label{thm:RCDstable}
  Let $\{(X_i,d_i,\m_i)\}_{i\in\bbN}$ be a sequence of $\RCD(K_i,N)$ spaces and $x_i^0\in\supp\m_i$. Assume $(X_i,d_i,\m_i,x_i^0)\xrightarrow{mGH}(X_{\infty},d_{\infty},\m_{\infty},x^0_{\infty})$ and assume $K_i\rightarrow K$. Then $(X_{\infty}, d_{\infty},\m_{\infty})$ is an $\RCD(K,N)$ space. Moreover, any sequence of $\RCD(K,N)$ spaces has a convergent subsequence in mGH sense. 
 \end{thm}
Let $(X,d,\m)$ be an $\RCD(K,N)$ space. Given $r>0$, $d_r:=r^{-1}d$ and $\m_r^x$ for $x\in \supp\m$ is defined by 
\begin{align}
 \m_r^x:=\left(\int_{B_r(x)}1-\frac{d(x,\cdot)}{r}\,d\m\right)^{-1}\m.\notag
\end{align}
Since $(X,d_r,\m_r^x)$ is an $\RCD(r^2K,N)$ space, combining Remark \ref{rem:domnmul} with Theorem \ref{thm:RCDstable}, we are able to find a convergent subsequence $\{(X,d_{r_i},\m^x_{r_i},x)\}_i$. Therefore we reach the following definition. 
 \begin{defn}[Tangent cone]
  For a given point $x\in \supp\m$, we define \emph{the tangent cone} at $x$ by 
  \begin{align}
   &\Tan(X,d,\m,x)\notag\\
   &:=\left\{(Y,d_Y,\m_Y,y)\;;\;(X,d_{r_i},\m^x_{r_i},x)\xrightarrow{mGH}(Y,d_Y,\m_Y,y)\text{ for a sequence }r_i\downarrow 0\right\}.\notag
  \end{align} 
  We often denote it by $\Tan(X,x)$ for short. 
 \end{defn}
The $\ell$-dimensional regular set $\calR_{\ell}$ is defined by 
\begin{align}
 \calR_{\ell}:=\left\{x\in X\;;\; \Tan(X,x)=\{(\R^{\ell},d_E,\ul{\calL^{\ell}},0)\}\right\},\notag
\end{align}
where $d_E$ is the standard Euclidean distance and $\ul{\calL^{\ell}}$ is the normalized Lebesgue measure on $\R^{\ell}$, this means, 
\begin{align}
 \int_{B_1(0)}1-\lv x\rv\,d\ul{\calL^{\ell}}(dx)=1.\notag
\end{align}
It is known that $\calR_{\ell}=\emptyset$ if $\ell>[N]$ for $\RCD(K,N)$ spaces.

 Brue and Semola proved the following result. 
 \begin{thm}[\cite{BSconst}]\label{thm:BS}
  Let $(X,d,\m)$ be an $\RCD(K,N)$ space. Then there exists an integer $n$ such that $\m(X\setminus \calR_n)=0$. 
 \end{thm}
 \begin{defn}[Essentially dimension]\label{def:essdim}
  We call the integer $n$ in Theorem \ref{thm:BS} \emph{the essential dimension}, and denote it by $\dime(X,d,\m)$ and $\calS:=X\setminus \calR_n$ the singular set in $X$. 
 \end{defn}
 \begin{rem}[Hausdorff dimension and Essential dimension]
  In \cite{Stmms2}, Sturm proved that the Hausdorff dimension is at most $N$ for $\RCD(K,N)$ spaces. By the behavior of measure on the regular set (see \cite{AHTweyl}), it is clear that $\dime\leq \mathsf{dim}_H$. The coincident of these two notion of dimension was open. However, recently Pan and Wei proved that there exists an $\RCD(K,N)$ space whose Hausdorff dimension is strictly larger than essential one (\cite{PW}).  
 \end{rem}
Theorem \ref{thm:BS} does not guarantee the non-existence of points belonging to another dimensional regular set(Non-existence of higher dimensional regular point is proven in \cite{K}). One dimensional case is much simpler than the situation for other dimension. 
 \begin{thm}[\cite{KL,CMisop,CMnew}] 
  Let $(X,d,\m)$ be an $\RCD(K,N)$ space. Then the following are all equivalent to each other; 
  \begin{enumerate}
   \item $\calR_1\neq\emptyset$, 
   \item $\dime(X,d,\m)=1$, 
   \item $(X,d)$ is isometric to either $\R$, $\R_{\geq 0}$, $[0,\ell]$, or $S^1(r)$ for $\ell>0$, $r>0$. 
  \end{enumerate}
  Moreover the reference measure $\m$ is equivalent to $\calH^1$(denote it by $\m\sim\calH^1$), which is of the form $\m=e^{-f}\calH^1$, and its density function $f$ is  $(K,N-1)$-convex. Hence $h:=e^{-f}$ satisfies (\ref{eq:hratiioest}). 
 \end{thm}
Recall we say that two measures $\sigma$ and $\tau$ are equivalent to each other if both $\sigma\ll\tau$ and $\tau\ll\sigma$ hold.  
 \begin{rem}\label{rem:contih}
  $(K,N)$-convexity of the density function is proven in \cite{KL}. The improvement version of the convexity is proven by \cite{CMonge,CMil}. The density function $h:=e^{-f}$ is continuous.
 \end{rem}
\subsection{Weyl's law on finite dimensional compact $\RCD$ spaces}
In this subsection, we always assume the metric measure space $(X,d,\m)$ is a compact $\RCD(K,N)$ space. As aforementioned before, all the spectrum of $-\Delta$ are eigenvalues, and $0=\lambda_0<\lambda_1\leq \lambda_2\leq \cdots\rightarrow \infty$ holds. We define the counting function $N(\lambda)$ by 
\begin{align}
 N(\lambda):=\#\left\{i\;;\;\lambda_i\leq \lambda\right\}.\notag
\end{align}
Ambrosio, Honda, and Tewodrose \cite{AHTweyl} proves the following result. 
 \begin{thm}[\cite{AHTweyl} cf. \cite{ZZweyl}]
  Let $(X,d,\m)$ be a compact $\RCD$ space with $\dime(X,d,\m)=n$. Then 
  \begin{align}
   \int_X\lim_{r\rightarrow+0}\frac{r}{\m(B_r(x))}\,\m(dx)=\lim_{r\rightarrow+0}\int_X\frac{r}{\m(B_r(x))}\m(dx)<\infty\notag
  \end{align}
  if and only if 
  \begin{align}
   \lim_{\lambda\rightarrow\infty}\frac{N(\lambda)}{\lambda^{n/2}}=\frac{\omega_n}{(2\pi)^n}\calH^n(\calR_n).\notag
  \end{align}
 \end{thm}
\begin{rem}
 Independently \cite{ZZweyl} also proves a similar result.  
\end{rem}
For $k\in\bbN$, define the subset $\calR_k^*$ by 
\begin{align}
 \calR_k^*:=\left\{x\in\calR_k\;;\;~^{\exists}\lim_{r\downarrow 0}\frac{\m(B_r(x))}{\omega_kr^k}\in(0,\infty)\right\}.\notag
\end{align}
It is known that $\m(\calR_k\setminus\calR_k^*)=0$ (see \cite{AHTweyl}). In order to prove the main result, we need the following results. 
 \begin{thm}[\cite{AHTweyl}]\label{thm:liminfineq}
  Let $(X,d,\m)$ be a compact $\RCD(K,N)$ space and $k=\dime(X,d,\m)$. Then we have 
  \begin{align}
   \liminf_{t\rightarrow +0}\left(t^{k/2}\sum_ie^{-\lambda_it}\right)\geq \frac{1}{(4\pi)^{k/2}}\calH^k(\calR_k^*)>0,\notag
  \end{align}
  where $\calH^k$ is the $k$-dimensional Hausdorff measure on $(X,d)$. 
 \end{thm}
The so-called \emph{Abelian theorem} is also important for our main result. 
 \begin{thm}[Abelian theorem cf.\cite{AHTweyl}]
  Let $\nu$ be a nonnegative and $\sigma$-finite Borel measure on $[0,\infty)$. Assume that there exist $\gamma\in [0,\infty)$ and $C\in [0,\infty)$ such that 
  \begin{align}
   \lim_{a\rightarrow \infty}\frac{\nu([0,a])}{a^{\gamma}}=C.\notag
  \end{align} 
  Then 
  \begin{align}
   \lim_{t\rightarrow+0}t^{\gamma}\int_{[0,\infty)}e^{-tx}\,d\nu(x)=C\Gamma(\gamma+1).\notag
  \end{align}
 \end{thm}
In the next section, we use the above theorem for $\nu=\sum_i\delta_{\lambda_i}$. Note that 
\begin{align}
 \int_{[0,\infty)}e^{-tx}\,d\nu(x)=\sum_ie^{-\lambda_it}\notag
\end{align} 
in this case. 
%
%
\section{Proof of the main theorem}
Now let $(X,d,\m)$ be a compact $\RCD(K,N)$ space with $\m(X)=1$. Without loss of generality, we may assume $K<0$ and $-K=N-1$ for simplicity. 
\medskip
\par The equivalence between (1) and (2) is proven in \cite{KL}. The implication (4) to (3), (3) to (5) are trivial. 
\medskip
\par \ul{(2) $\Rightarrow$ (4)}: When $X=S^1(r)$, then the density function $h$ for $\m=h\calH^1$ is continuous(see Remark \ref{rem:contih}). Since $\m\sim \calH^1$ and the continuity of $h$, $h$ never vanish. Put $c:=\inf h>0$ and $C:=\sup h<\infty$. Since 
\begin{align}
 0<\frac{r}{\m(B_r(x))}h(x)=\frac{rh(x)}{\int_{B_r(x)}h(t)\,dt}\leq \frac{rh(x)}{2rc}\leq \frac{C}{2c}<\infty\notag
\end{align} 
holds, we have 
\begin{align}
 \lim_{r\rightarrow+0}\int_X\frac{r}{\m(B_r(x))}\,\m(dx)=\int_X\lim_{r\rightarrow+0}\frac{r}{\m(B_r(x))}\,\m(dx)<\infty\label{eq:S1dom}
\end{align}
by the dominated convergence theorem. 
\smallskip
\par 
Let us consider the case for $X=[0,\ell]$. Without loss of generality, we may assume $\ell=\pi$. When the density function $h$ for $\m=h\calH^1$ has a positive minimum, then a similar argument as $X=S^1$ implies the consequence. So, we assume $h(0)=h(\pi)=0$. Take a small positive number $r>0$ and fix it. For $x\in [0,r)$, we have 
\begin{align}
 h(t(x+r))^{1/N-1}&\geq \sigma^{(1-t)}_{K,N-1}(x+r)h(0)^{1/N-1}+\sigma^{(t)}_{K,N-1}(x+r)h(x+r)^{1/N-1}\notag\\
 &=\frac{\sinh (t(x+r))}{\sinh(x+r)}h(x+r)^{1/N-1}.\notag
\end{align}
Thus we get 
\begin{align}
 &\m(B_r(x))=\int_0^{x+r}h(y)\,dy=(x+r)\int_0^1h(t(x+r))\,dt\notag\\
 &\geq (x+r)\int_0^1\frac{\sinh^{N-1}(t(x+r))}{\sinh^{N-1}(x+r)}h(x+r)\,dt.\notag
\end{align}
Therefore 
\begin{align}
 0<\frac{r}{\m(B_r(x))}\leq \frac{r\sinh^{N-1}(x+r)}{(x+r)h(x+r)}\left(\int_0^1\sinh^{N-1}(s(x+r))\,ds\right)^{-1}.\notag
\end{align}
We have 
\begin{align}
 &\frac{rh(x)}{\m(B_r(x))}\notag\\
 &\leq \frac{rh(x)\sinh^{N-1}(x+r)}{(x+r)h(x+r)}\left(\int_0^1\sinh^{N-1}(s(x+r))\,ds\right)^{-1}\notag\\
 &=\frac{r}{x+r}\frac{h(x)}{h(x+r)}\frac{\sinh^{N-1}(x+r)}{\int_0^1\sinh^{N-1}(s(x+r))\,ds}.\notag
\end{align}
By using the Taylor expansion for $\sinh$, we have $\sinh^{N-1}(s(x+r))\geq s^{N-1}(x+r)^{N-1}$ and $\sinh^{N-1}z\leq 2^{N-1}z^{N-1}$ for sufficiently small $z>0$, further, applying $(\ref{eq:hratiioest})$ to $h$ on $(0,x+2r)$, we obtain 
\begin{align}
 \frac{h(x)}{h(x+r)}\leq \frac{\sinh^{N-1}((x+2r)-x)}{\sinh^{N-1}((x+2r)-(x+r))}=\frac{\sinh^{N-1}(2r)}{\sinh^{N-1}(r)}\leq\frac{(4r)^{N-1}}{r^{N-1}}=4^{N-1}.\notag
\end{align}
Also we get 
\begin{align}
 &\sinh^{N-1}(x+r)\left(\int_0^1\sinh^{N-1}(s(x+r))\,ds\right)^{-1}\notag\\
 &\leq 2^{N-1}(x+r)^{N-1}\left(\frac{(x+r)^{N-1}}{N}\right)^{-1}\notag\\
 &\leq N2^{N-1}.\notag
\end{align}
Finally we have 
\begin{align}
 \frac{rh(x)}{\m(B_r(x))}\leq N8^{N-1}.\notag
\end{align}

On the other hand, for $x\in (r,\pi/2)$, we have 
\begin{align}
 &\m(B_r(x))=\int_{x-r}^{x+r}h(y)\,dy=\int_0^1h\left((1-t)(x-r)+t(x+r)\right)\cdot 2r\,dt\notag\\
 &\geq 2r\int_0^1\left(\sigma^{(1-t)}_{K,N-1}(2r)h(x-r)^{1/N-1}+\sigma^{(t)}_{K,N-1}(2r)h(x+r)^{1/N-1}\right)^{N-1}\,dt\notag\\
 &\geq 2r\int_0^1\sigma^{(t)}_{K,N-1}(2r)^{N-1}h(x+r)\,dt\notag\\
 &=\frac{2rh(x+r)}{\sinh^{N-1}\left(2r\right)}\int_0^1\sinh^{N-1}(2tr)\,dt\notag\\
 &\geq \frac{2rh(x+r)}{(4r)^{N-1}}\int_0^1(2tr)^{N-1}\,dt\notag\\
 &\geq \frac{2rh(x+r)}{N2^{N-1}}.\notag
\end{align}
Applying (\ref{eq:hratiioest}) to $h$ on $(x-r,x+2r)$, we get the estimate 
\begin{align}
 \frac{h(x)}{h(x+r)}\leq \left(\frac{\sinh((x+2r)-x)}{\sinh((x+2r)-(x+r))}\right)^{N-1}=\frac{\sinh^{N-1}(2r)}{\sinh^{N-1}(r)}\leq \frac{(4r)^{N-1}}{r^{N-1}}=4^{N-1}.\notag
\end{align}
Then 
\begin{align}
 \frac{rh(x)}{\m(B_r(x))}&\leq \frac{rh(x)N2^{N-1}}{2rh(x+r)}=\frac{N2^{N-1}}{2}\cdot \frac{h(x)}{h(x-r)}\notag\\
 &\leq \frac{N2^{N-1}}{2}\cdot 4^{N-1}\leq N8^{N-1}.\notag
\end{align}
The upper bound $N8^{N-1}$ depends on neither $x\in(0,\pi/2)$ nor $r>0$. We apply the same argument for near $x=\pi$. 
Then by the dominated convergence theorem, 
\begin{align}
 \lim_{r\downarrow 0}\int_X\frac{r}{\m(B_r(x))}\,\m(dx)=\int_X\lim_{r\downarrow0}\frac{r}{\m(B_r(x))}\,\m(dx)<\infty.\label{eq:elldom}
\end{align} 
Both cases, $X=S^1(r)$, $[0,\ell]$, we have (4) by combining (\ref{eq:S1dom}) and (\ref{eq:elldom}) with Theorem 4.3 in \cite{AHTweyl}.

\medskip
\par\ul{(5) $\Rightarrow$ (1)}: By combining the assumption and Abelian theorem, we have 
\begin{align}
 \lim_{t\rightarrow +0}t^{\frac{1+\alpha}{2}}\sum_ie^{-\lambda_it}=0.\notag
\end{align}
Let $k=\dime(X,d,\m)$. Then by Theorem \ref{thm:liminfineq}, we obtain 
\begin{align}
 \lim_{t\rightarrow+0}t^{\frac{1+\alpha}{2}}\sum_ie^{-\lambda_it}&=0<\frac{1}{(4\pi)^{k/2}}\calH^k(\calR_k^*)\notag\\
 &\leq \liminf_{t\rightarrow+0}t^{\frac{k}{2}}\sum_ie^{-\lambda_it}.\notag
\end{align}
This implies $1+\alpha> k$. Since $\alpha\leq 1$ and $k$ is an integer, $k$ has to be 1.

\section*{Acknowledgement}
The authors would like to thank Professor Shouhei Honda for telling us the problem and his helpful comments and fruitful discussion. Y.K. is partly supported by JSPS KAKENHI Grant Numbers JP18K13412 and JP22K03291.

\end{document}